\documentclass[a4,12pt,reqno,fleqn]{amsart}
\usepackage{graphicx}
\usepackage{amscd}
\usepackage{xfrac}
\usepackage{mathtools}
\usepackage[english]{babel}
\usepackage[autostyle]{csquotes}
\emergencystretch=\maxdimen
\allowdisplaybreaks[1]
\usepackage{array}
\usepackage{xcolor}
\usepackage[all]{xy}
\makeatletter
\@namedef{subjclassname@2020}{\textup{2020} Mathematics Subject Classification}
\makeatother


\usepackage{amsmath}
\usepackage{amsthm}
\usepackage{amsfonts,amssymb,hyperref,cleveref,enumitem,mathrsfs,setspace,verbatim}
\usepackage[margin=3.5cm]{geometry}

\DeclareMathAlphabet{\mathpzc}{OT1}{pzc}{m}{it}
\DeclareMathOperator{\sign}{sign}
\newtheorem{thm}{Theorem}[section]
\newtheorem{cor}[thm]{Corollary}

\newtheorem{propn}[thm]{Proposition}
\theoremstyle{definition}
\newtheorem{defn}[thm]{Definition}

\newtheorem{example}[thm]{Example}

\newcommand{\bj}{\perp_{BJ}}
\newcommand{\nbj}{\not\perp_{BJ}}

\newcommand{\ot}{\otimes}
\newcommand{\Pg}{\mathcal{P}_G}
\author[Mohit and R.~Jain]{Mohit and Ranjana Jain}

\address{Mohit, Department of  Mathematics, University of Delhi, Delhi}
\email{mohitdhandamaths@gmail.com}

\address{Ranjana Jain, Department of Mathematics, University of Delhi, Delhi}
\email{rjain@maths.du.ac.in}
\thanks{Research of the first named author is supported by Savitribai Jyotirao Phule Single Girl Child Fellowship vide F.No. 82-7/2022(SA-III) and the second named author is supported by Faculty Research Programme Grant by IoE, University of Delhi vide Ref.No./IoE/2023-24/12/FRP}

\begin{document}
	\title[BJ-orthogonality in tensor product of Lebesgue-Bochner spaces]{Birkhoff-James orthogonality in certain tensor products of Banach spaces II}
	\maketitle
	
		\textbf{Abstract:} In this article, we discuss the relationship between Birkhoff-James orthogonality of elementary tensors in the space $L^{p}(\mu)\otimes^{\Delta_{p}}X,\; (1\leq p<\infty)$ with the individual elements in their respective spaces, where $X$ is a Banach space whose norm is Fr$\acute{e}chet$ differentiable and $\Delta_{p}$ is the natural norm induced by $L^{p}(\mu,X)$. In order to study the said relationship, we first provide some characterizations of Birkhoff-James orthogonality of elements in the Lebesgue-Bochner space $L^{p}(\mu,X)$.
	\section{Introduction}

The notion of Birkhoff-James orthogonality (in short, BJ-orthogonality) has been studied extensively in the last few decades  in the category of normed spaces and more recently it has also attracted researchers from the areas of Banach algebras and operator algebras. 
One of the major reason behind the study of the notion of BJ-orthogonality is that it has various applications in the geometry of Banach spaces, see \cite{birk,james}.  Recall that, for a normed space $X$ over a field $\mathbb{F}\ (\mathbb{R}$ or $\mathbb{C}$) and $x,\;y\in X,\;x$ is said to be BJ-orthogonal to $y$, denoted as $x \bj y$  if    
\begin{center}
	$	\|x+\lambda y\|\geq \|x\|,$ for all $\lambda\in\mathbb{F}.$
\end{center}  
	Interestingly, there exists a very natural relationship between the notion of BJ-orthogonality and the (more classical and much thoroughly studied) notion of  best approximation elements in normed spaces.  The problem of best approximation lies in finding, for a given $x$ in a normed space $X$, an element $g_{0}$ in a subset (mostly a subspace) $G$  of $X$ such that  
$$
\|x-g_{0}\| = \inf\{ \|x-g\| :  g\in G\}.
$$
 Such an element  $g_{0}$ is called a point of best approximation of $x$ out of $G$.  We shall represent  the set of all best approximants of $x$ in $G$ by $\mathcal{P}_{G}(x)$. It is rather straight forward to see that for any two elements $x,y$ in a normed space $X$,  $x \bj y$ if and only if $0 \in \Pg(x)$, $G$ being the subspace spanned by $y$. Needless to mention, a good deal of research has been done from both perspectives.

Since the seminal work of Grothendieck on tensor products of topological vectors spaces \cite{gro}, the theories of tensor products of Banach spaces and Banach algebras (in particular) have proved to be indispensible in the proper  understanding of these categories. Our motivation to study BJ-orthogonality arose from the natural question of studying tensor product spaces from the perspective of BJ-orthogonality.  In this direction, in \cite{jain}, we initiated the analysis of  BJ-orthogonality of elementary tensors in various tensor product spaces in terms of  the BJ-orthogonality of the  individual elements in their respective spaces. One of the main results exhibited that for real Banach spaces $X$ and $Y$, $x_{1}\otimes y_{1}\perp_{BJ}x_{2}\otimes y_{2}$ in $X\otimes^{\lambda}Y$ if and only if $x_{1}\perp_{BJ}x_{2}$ or $y_{1}\perp_{BJ}y_{2},$ where $\lambda$ is the injective tensor product \cite[Theorem 3.7]{jain}. However this equivalence fails in $C_{\mathbb{C}}(K_{1})\otimes^{\alpha}C_{\mathbb{C}}(K_{2})$, for any reasonable cross norm $\| \cdot \|_{\alpha}$ \cite[Example 3.2]{jain}. It was also established that this equivalence is true in $ L^p(\mu) \otimes^{\Delta_{p}} L^p(\nu)$, $1 < p<\infty$, where $\mu$, $\nu$ are positive measures and $\Delta_{p}$ is defined in \Cref{norm} (see, \cite[Theorem 3.9]{jain}). 
 
Having made these observations, it was then quite natural to ask whether similar results hold in the tensor product spaces $L^1(\mu)\otimes^{\Delta_p} X$ for $1\leq p < \infty$. To our somewhat mixed satisfaction, in this  article, we establish that a similar result holds in  $L^{p}(\mu)\otimes^{\Delta_{p}} X$ as well, where $X$ is a Banach space whose norm is Fr$\acute{e}$chet differentiable, $\mu$ is a positive measure and $1 < p < \infty$.  When $p=1$,  quite surprisingly, we observe through concrete examples that such a result is not true in the tensor product spaces $L^{1}(\mu)\otimes^{\gamma} L^{1}(\nu),$ and  $L^{1}(\mu)\otimes^{\gamma} X$, where $\otimes^{\gamma}$ is the Banach space projective tensor product. 
		
		
		Our approach relies on the well known isometric identification  between $L^p(\mu) \otimes^{\Delta_p} X$ and the Lebesgue-Bochner space of $p$-integrable functions  $L^p(\mu, X)$ (see, \cite[$\S$7.2]{floret}) and an appropriate exploitation of the above mentioned relationship between the notions of BJ-orthogonality and best approximation elements. More precisely,   motivated by some classical results of Kripke-Rivlin and Singer regarding best approximations (as mentioned in the following two paragraphs) in $L^p$-spaces, we first obtain certain characterizations of best approximation elements  for points in $L^p(\mu, X)$  and then we exploit those to derive the above mentioned results in the context of tensor products. 
		
		In 1965, Kripke and Rivlin \cite{kripke}  characterized the elements of $\Pg(f)$ for $f\in L^{1}(\mu)\setminus\overline{G}$,  
		 for a positive measure space $(S,\mu)$ and a subspace $G$ of $L^{1}(\mu)$. In particular, he proved that $g_0\in \mathcal{P}_G(f)$ if and only if
		 \begin{equation}\label{J} 
		\hspace{1.5cm} \bigg|\int\limits\limits\limits_{S}{g(s) \ \overline{\sign}(f(s)-g_0(s))} \,ds\bigg|\leq \int\limits\limits\limits_{Z(f-g_o)}{|g(s)| \,ds}
		\end{equation}
		for all $g\in G$, here $Z(f)$ denotes the zero set $\{s\in S:f(s) =0\}$ of $f$ and $\sign (a)=	a/|a| $ for $0 \neq a \in \mathbb{C},\; \sign(0)=0$.
		
		Later, in 1970, Singer \cite{singer} gave a different proof of this characterization, and also provided its analogue for elements in $L^{p}(\mu)$, $1<p< \infty$.  More precisely, he proved that if $(S,\mu)$ is a positive measure space and $G$ is a subspace of $L^{p}(\mu)$, then for $f\in L^{p}(\mu)\setminus \overline{G},\; g_0 \in \mathcal{P}_G(f)$ if and only if
		$$\int\limits\limits\limits_{S}{g(s)|f(s)-g_0(s)|^{p-1}\overline{sign}(f(s)-g_0(s))} \,ds=0$$
		for all $g\in G$ \cite[Theorem 1.11]{singer}. A natural question arises that what can be said about the approximation for the vector valued integrable functions. In 1989, Smirnov \cite{smirnov} characterized the elements of $\Pg(f)$ for $f\in C([0,1], X)\setminus U$, for a real smooth Banach space $X$ and a convex subset $U$ of $C([0,1], X)$(equipped with the integral-norm).

		In Section 3, for any positive measure space $(S, \Sigma, \mu$) and a  Banach space $X$  belonging to a relatively large class, we provide natural analogues of the characterizations of best approximation elements by Kripke-Rivlin and Singer in $L^{p}(\mu, X), 1\leq p < \infty $, using entirely different techniques.  
		Before using these characterizations in the tensor product context, for $p>1$, we derive a new and quite elegant proof of Light's theorem \cite[Corollary 2]{light} for Banach spaces with Fr$\acute{e}$chet differentiable norm. More precisely, we prove that  if $f\in L^{p}(\mu,X)$ $(1< p<\infty)$ and $Y$ is a closed subspace of $X$, where $X$ is a Banach space whose norm is Fr$\acute{e}$chet differentiable and  $\mu$ is a finite complete positive measure, then $f\bj L^{p}(\mu,Y)$ if and only if $f(s)\bj Y$ for a.e.\enquote{$s$}.  We must mention that the  proof by  Light  was based on some distance formula. 
		
Finally, in Section 4, we present our main results related to the BJ-orthogonality of tensor products. 


\section{Preliminaries}

	We first collect some basic definitions and results which we need for our investigation.
 For a measure space $(S, \Sigma,\mu)$ ($\mu$ being a positive measure) and a Banach space $X$ (real or complex), {\it Lebesgue-Bochner space} is defined as
 $$  \hspace*{-0.3cm} L^{p}(\mu, X)= \{f:S\rightarrow X|\ f \ \text{is strongly measurable and} \ \;\int\limits\limits\limits_{S}\|f(s)\|^p\,ds <\infty\},
  $$
$1\leq p<\infty,$ where almost everywhere equal functions are identified.

 \begin{defn}\label{norm}
 	For $u\in L^{p}(\mu)\otimes X$, define the  $\Delta_{p}$-norm as 
$$ 		||u||_{\Delta_{p}}=||\phi(u)||_{L^{p}(\mu, X)}, $$ 
 where $\phi:L^{p}(\mu)\otimes X\rightarrow L^{p}(\mu, X)$ given by $\phi(f\otimes x)=f(\cdot)x$ is an injective map. We denote  the completion of  $L^{p}(\mu)\otimes X$ with respect to the $\Delta_{p}$-norm by $L^{p}(\mu)\otimes^{\Delta_{p}} X$.
\end{defn}

 For $p=1$, $\Delta_{p}$-norm coincides with the {\it Banach space projective tensor norm} given by
 $$ \|u\|_{\gamma} = \inf \left\{\sum_{i=1}^{n} \|f_i\| \|x_i \|; \ u =\sum_{i=1}^{n} f_i \ot x_i \right\}, \ \ u \in L^{1}(\mu)\otimes X .$$
It is well known that the spaces $L^{p}(\mu)\otimes^{\Delta_{p}}X$ and $L^1(\mu)\otimes^{\gamma}L^1(\nu)$ are isometrically isomorphic to $L^{p}(\mu,X)$ and $L^1(\mu\times\nu)$, respectively see \cite[$\S$ 7.1, 7.2]{floret}.

	For $0\neq x\in X$, a {\it support map} $F_{x}$ at $x$ is a bounded linear functional on $X$ of norm one satisfying $F_{x}(x)=\|x\|.$  An element $x$ is said to be {\it smooth} if $F_x$ is unique, and if every non-zero element of $X$ is smooth then we call $X$ to be {\it smooth}.
	
	The norm function $\| \cdot \|$ on $X$ is said to be
	\begin{itemize}
		\item  {\it G$\hat{a}$teaux differentiable} at  $0\neq x\in X$ if $$\lim_{\alpha \to 0}\frac{\|x+\alpha y\|-\|x\|}{\alpha}$$ exists for every $y\in X$.	It is well known that an element $x$ is smooth if and only if the norm is G$\hat{a}$teaux differentiable at  $x$, and in this case, the G$\hat{a}$teaux derivative at $x$, takes the value $Re(F_{x}(y))$ in the direction of $y\in X$. 
		
		\item {\it Fr$\acute{e}$chet differentiable} at $0\neq x\in X$ if there exists an $f \in X^*$  satisfying $$\lim_{h \to 0}\frac{\big|\|x+h\|-\|x\|-f(h)\big|}{\|h\|}=0.$$ The norm function is Fr$\acute{e}$chet differentiable if it is Fr$\acute{e}$chet differentiable at every non-zero point.
	\end{itemize}
	
It is known  that Fr$\acute{e}$chet differentiability implies G$\hat{a}$teaux differentiablity, and thus  the space $X$ is smooth,  if the norm function on $X$ is Fr$\acute{e}$chet differentiable.


Throughout the article, $(S,\mu)$ denotes a positive measure space unless specified.

	\section{Best Approximation in $L^p(\mu,X)$}
	 We first suitably characterize the elements of best approximations  in a subspace of $L^{1}(\mu,X)$. One can easily observe that for $X = \mathbb{F}$, by Riesz representation theorem, support map $F_{z_0}$ is given by $F_{z_0}(z)=\left\langle z,\frac{z_0}{|z_0|}\right\rangle_X$=$z\cdot\overline{\sign}(z_0)$. Thus, the Kripke-Rivlin's characterization given in \ref{J} can be reformulated as $ g_0\in \mathcal{P}_G(f)$ if and only if
	 \begin{equation*}
	 \hspace{1.5cm}	\bigg|\int\limits\limits\limits_{Z(f-g_0)^c}{F_{f(s)-g_0(s)}( g(s)) \,ds}\bigg|\leq\int\limits_{Z(f-g_0)}{|g(s)| \,ds}
	 \end{equation*}
	 for every $g\in G.$ It is worth mentioning that we get a similar characterization for vector-valued integrable functions for a real Banach space when the norm on $X$ is Fr$\acute{e}$chet differentiable. Further, when $X$ is a complex  Banach 
space we have a slightly weaker charaterization, but for a large class of spaces, namely smooth Banach spaces.  We would deploy few techniques of Smirnov \cite{smirnov} to prove the necessary part of the following result.
	 \begin{thm}\label{L1chr}
	 	Let $X$ be a complex smooth Banach space and $G$ be a subspace of $L^{1}(\mu, X)$. For $f\in\;L^{1}(\mu,X)\setminus \overline{G},\; g_0\in \mathcal{P}_G(f)$ if and only if
	 	\begin{equation}\label{1}
	 		\bigg|\int\limits\limits\limits_{Z(f-g_0)^c}{Re(F_{f(s)-g_0(s)}( g(s))) \,ds}\bigg|\leq\int\limits\limits\limits_{Z(f-g_0)}{\|g(s)\| \,ds}
	 	\end{equation}
	 	for all $g\in G.$
	 \end{thm}
	 \begin{proof}
	 	First suppose that $g_0\in \mathcal{P}_G(f)$ and let $g$ be any arbitrary element of $ G$. Then for $\alpha \in \mathbb{C}$, $\|f-g_0+\alpha(g_0-g)\|_{1}\geq \|f-g_0\|_{1}$, that is,
	 	\begin{equation}\label{2c}
	 		\int\limits\limits\limits_{S} \|f(s)-g_0(s)+\alpha(g_0(s)-g(s))\|_{X}-\|f(s)-g_0(s)\|_{X} \, ds \geq 0. \end{equation}	 
	 	For each $n\in\mathbb{N}$, consider a measurable function $h_n:S\rightarrow \mathbb{R}$ given by $$h_n(s)= n\big(\big\|f(s)-g_0(s)+\frac{1}{n}(g_0(s)-g(s))\big\|_{X}-\big\|f(s)-g_0(s)\big\|_{X}\big),  \ s \in S.$$	 
	 	By the triangle inequality we have $|h_n(s)|\leq\|g_0(s)-g(s)\|$ for all $s\in S$ and $n\in\mathbb{N}$. 
	 	Thus, by Lebesgue dominated convergence theorem and by (\ref{2c}), we have $   \int\limits\limits\limits_{S}\lim_{n \to \infty}{h_n(s) \, ds}\geq 0. $	Now, for  $s\in S \setminus Z(f-g_0)$, G$\hat{a}$teaux differentiability of the norm function gives
	 	$$\lim_{n \to \infty}{h_n(s)}=Re(F_{f(s)-g_0(s)}(g_0(s)-g(s))).$$
	 	Thus, we have
	 	$$\int\limits\limits\limits_{Z(f-g_0)^c}{Re(F_{f(s)-g_0(s)}(g_0(s)-g(s)) )\,ds}+\int\limits\limits\limits_{Z(f-g_0)}{\|g_0(s)-g(s)\|}\,ds\geq 0.$$
	 	for all $g\in G$.
	 	Since $G$ is a subspace, replacing $g$ by $g_0-g$ and $g_0+g$ and using the fact that real part of support map is real linear 
	 	we have
	 	$$\bigg|\int\limits\limits\limits_{Z(f-g_0)^c}{Re(F_{f(s)-g_0(s)}(g(s))) \,ds}\bigg|\leq\int\limits\limits\limits_{Z(f-g_0)}{\|g(s)\| \,ds}$$
	 	for all $g\in G$.
	 	
	 	Conversely, assume that the inequality holds. In order to prove $g_0\in \mathcal{P}_G(f)$, equivalently, $f-g_0\perp_{BJ}G$, by \cite[Proposition 1.5]{keckic}, it is sufficient to prove that $\underset{\phi\in[0,2\pi)}{\it{inf}}D_{\phi,\;f-g_0}(g)\geq 0$ for all $g\in\;G$, where $D_{\phi,\;f-g_0}(g)$ is the  $\phi$-G$\hat{a}$teaux derivative of the norm function at point $f-g_0$ in the direction of $g$.  For this, let $g\in G$ be any arbitrary element and $\phi\in [0,2\pi).$ Then
	 	\begin{equation*}
	 		\begin{split}
	 			\hspace*{-1.3cm} D_{\phi,\;f-g_0}(g) & =  \lim_{\alpha \to 0^{+}}\frac{\|f-g_0+\alpha e^{i\phi} g\|_{1}-\|f-g_0\|_{1}}{\alpha}\\
	 			& \hspace*{-0.7cm} =  \lim_{\alpha \to 0^{+}}\int\limits\limits\limits_{S}\frac{1}{\alpha}(\|f(s)-g_0(s)+\alpha e^{i\phi} g(s)\|_{X}-\|f(s)-g_0(s)\|_{X}) \,ds\\	  
	 			& \hspace*{-0.7cm} = \int\limits\limits\limits_{S}\lim_{\alpha \to 0^{+}}\frac{1}{\alpha}(\|f(s)-g_0(s)+\alpha g_{1}(s)\|_{X}-\|f(s)-g_0(s)\|_{X})\,ds
	 		\end{split}
	 	\end{equation*}
	 	using Lebesgue dominated convergence theorem, where $g_{1}=e^{i\phi}g\in G$ as $X$ is a complex Banach space.
	 	Using (\ref{1}) for $g_{1}$, we have 
	 	$$\int\limits\limits\limits_{Z(f-g_0)^c}{Re(F_{f(s)-g_0(s)}(g_{1}(s)) )\,ds}+\int\limits\limits\limits_{Z(f-g_0)}{\|g_{1}(s)\|}\,ds\geq 0.$$
	 	Since the norm function on $X$ is G$\hat{a}$teaux differentiable, the above inequality reduces to
	 	
	 	$$\int\limits\limits\limits_{S}\lim_{\alpha \to 0^{+}}{\frac{\|f(s)-g_0(s)+\alpha g_{1}(s)\|_{X}-\|f(s)-g_0(s)\|_{X}}{\alpha}\,ds}\geq 0,$$
	 	which completes the proof.
	 \end{proof}
 Next, we derive a similar characterization when $X$ is a real Banach space. To do this, we first prove a result which is motivated from \cite[Lemma 3]{giles}, with a slightly different notion of support map. Also, it is worth mentioning that for the real Banach spaces, this result was first proved by Cudia \cite[Corollary 4.11]{cudia}. However, we present a much simpler proof for any Banach space.  
	 \begin{propn}\label{continuity}
	 	Let $X$ be a Banach space whose norm is Fr$\acute{e}$chet differentiable. Then  the map $T:X\rightarrow X^{*}$ defined by $T(x)=F_{x}$ for $0\neq x$ and $T(0)=0$ is continuous on $X\setminus\{0\}$, where $X$ and $X^{*}$  are equipped with norm topologies. 
	 \end{propn}
	 \begin{proof}
	 	 Since $F_{\alpha x}=F_{x}$ for $\alpha >0$, it is enough to prove that $T$ is continuous on the unit sphere $S_{1}$. Consider an arbitrary element $x\in S_{1}$. We first claim that $T$ is continuous at $x$, when $X^*$ is equipped with the weak$^*$-topology.\\
	 	 If not, then there exist a weak$^*$-neighbourhood, say $V$, of $T(x)$ and a sequence $\{y_{n}\}$ in $X$ such that $\|x-y_{n}\|\leq \frac{1}{n}$ and $T(y_{n})\notin V.$ Since $\underset{n}\lim\|y_{n}\|=\|x\|=1$, there exists $m\in\mathbb{N}$ such that $\|y_n\|\neq0$ for all $n\geq m$. thus without loss of generality, we consider the sequence $\{y_n\}$ of non-zero terms. By Banach-Alaoglu theorem, the closed unit ball $B^{*}$ of $X^*$ is weak$^*$-compact and thus the sequence $\{F_{y_{n}}\}$ in $B^{*}$ has a convergent subnet say $\{F_{y_{\theta(k)}}\}$, converging to $f\in B^*$.  Let $c$ be a bound of the sequence  $\{\frac{1}{\|y_{n}\|}\}$. We claim that $f=F_{x}$. To see this, consider
	 	\begin{align*}
	 		|f(x)-1|&=\bigg|f(x)-F_{y_{\theta(k)}}\bigg(\frac{y_{\theta(k)}}{\|y_{\theta(k)}\|}\bigg)\bigg|\\
	 		&\leq|f(x)-F_{y_{\theta(k)}}(x)|+\bigg|F_{y_{\theta(k)}}(x)-F_{y_{\theta(k)}}\bigg(\frac{y_{\theta(k)}}{\|y_{\theta(k)}\|}\bigg)\bigg|\\
	 		&\leq|f(x)-F_{y_{\theta(k)}}(x)|+\frac{1}{\|y_{\theta(k)}\|}\big\|(x\|y_{\theta(k)}\|-y_{\theta(k)}\big)\|\\
	 		&\leq|f(x)-F_{y_{\theta(k)}}(x)|+c\big\|(x\|y_{\theta(k)}\|-y_{\theta(k)})\big\|,
	 	\end{align*} 
	 	Since $f$ is a weak$^*$-limit of $\{F_{y_{\theta(k)}}\}$, the subnet $\{F_{y_{\theta(k)}}(x)\}$ converges to $f(x)$. Also, the subnet $\{x\|y_{\theta(k)}\|-y_{\theta(k)}\}$ converges to zero. Thus, $f(x)=1=\|x\|$, and smoothness of $X$ gives $f=F_{x}$, that is, $F_{x}$ is a weak$^*$-cluster point of $\{F_{y_{n}}\}$. Thus, $V$ contains some points of the sequence $\{F_{y_{n}}\}$ which is a contradiction.\\
	 	Now, let $\{x_{n}\}$ be a sequence in $X$ converging to $x$. Then $T(x_n)\rightarrow T(x)$ in weak$^*$-topology of $X^{*}$. Thus, we have a sequence $\{F_{x_{n}}\}$ satisfying $\underset{n}\lim F_{x_{n}}(x)=\|x\|$. By \cite[Lemma 4]{giles}, $\{F_{x_{n}}\}$ is norm convergent to $F_{x}$. Hence $T$ is continuous at $x$, which completes the proof.
	 \end{proof}
  \begin{thm}\label{L1chreal}
	Let $X$ be a real Banach space whose norm is Fr$\acute{e}$chet differentiable and $G$ be a subspace of $L^{1}(\mu, X)$. For $f\in\;L^{1}(\mu,X)\setminus \overline{G},\; g_0\in \mathcal{P}_G(f)$ if and only if
	\begin{equation}\label{2}
		\bigg|\int\limits\limits\limits_{Z(f-g_0)^c}{F_{f(s)-g_0(s)}( g(s)) \,ds}\bigg|\leq\int\limits\limits\limits_{Z(f-g_0)}{\|g(s)\| \,ds}
	\end{equation}
	for all $g\in G.$
\end{thm}
\begin{proof}
	Since $X$ is smooth, the norm being Fr$\acute{e}$chet differentiable, proof of the necessary part follows on the same lines of \Cref{L1chr}.\\
For the converse, consider an arbitrary element $g$ of $G$.\\
	Case(i): If $\int\limits_{Z(f-g_{0})}\|g(s)\|\,ds=0$, then define $\phi:S\rightarrow X^{*}$ as
	\begin{equation*}
		\phi(s)=
		\begin{cases}
			F_{f(s)-g_{0}(s)} & \text{if}\;s\in {Z(f-g_{0})}^c,\\
			0 & \text{if}\;s\in Z(f-g_{0}).
		\end{cases}
	\end{equation*}
	 We first claim that $\phi\in L^{\infty}(\mu,X^{*})$. Since norm of the support map is one and $f\neq g_0$, it is sufficient to show that $\phi$ is strongly measurable.
	Since $f-g_{0}$ is strongly measurable, there exists a sequence $\{\psi_{n}\}$ of simple measurable functions such that $\psi_{n}(s)\rightarrow f(s)-g_{0}(s)$ for a.e. \enquote{$s$}. Let $A_n=\{s\in S: \psi_n(s)\neq0\}$ and consider a sequence $\{\phi_n\}$ of simple measurable functions defined as
	\begin{equation*}
		\phi_{n}(s)=
		\begin{cases}
			F_{\psi_{n}(s)} & \text{if}\;s\in {Z(f-g_{0})}^c\cap A_n,\\
			0, & elsewhere.
		\end{cases}
	\end{equation*}
	 Let $s\in S$ for which $\underset{n}\lim\psi_{n}(s)=f(s)-g_{0}(s)$.
	 If $s\in {Z(f-g_{0})}^c$, then there exists $n_0\in \mathbb{N}$ such that $s\in A_{n}$ for all $n\geq n_0$. Thus, by \Cref{continuity}, the sequence $\{F_{\psi_{n}(s)}\}_{n\geq n_0}$ converges to $F_{f(s)-g_{0}(s)}$  and hence the sequence $\{\phi_n(s)\}_{n\geq n_0}$ converges to $\phi(s)$.
	 If $s\in Z(f-g_{0})$, then $\phi_n(s)=0=\phi(s)$ for all $n$. Thus, in both the cases, the function $\phi$ is a.e. limit of a sequence of simple measurable functions and hence $\phi$ is strongly measurable function. Now, by (\ref{2}), we have $\int_{S}{\phi(s)(g(s))}\,ds=0$. Thus,
	 \begin{align*}
	 	\|f-g_0\|_1&=\int\limits_{S}{\|f(s)-g_0(s)\|}\,ds\\
	 	&=\int_{S}{\phi(s)(f(s)-g_0(s))}\,ds\\
	 	&\leq\int_{S}{|\phi(s)(f(s)-g_0(s)-g(s))|}\,ds\\
	 	&\leq\int_{S}{\|\phi(s)\|\|f(s)-g_0(s)-g(s)\|}\,ds\\
	 	&\leq\|\phi\|_{L^{\infty}(\mu,X^{*})}\int_{S}{\|f(s)-g_0(s)-g(s)\|}\,ds\\
	 	&=\|f-g_0-g\|_1.
	 \end{align*}
	Case(ii): If $\int\limits_{Z(f-g_{0})}\|g(t)\|\,dt\neq0$, set c= $\frac{-\big|\int\limits_{{Z(f-g_{0})}^c}F_{f(t)-g_{0}(t)}(g(t))\,dt\big|}{\int\limits_{Z(f-g_{0})}\|g(t)\|\,dt}$ so that by (\ref{2}), $|c|\leq1$. Define $\phi:S\rightarrow X^{*}$ as
	\begin{equation*}
		\phi(s)=
		\begin{cases}
			F_{f(s)-g_{0}(s)} & \text{if}\;s\in  Z(f-g_{0})^c,\\
			cF_{g(s)} & \text{if}\;s\in Z(f-g_{0})\cap {Z(g)}^c,\\
			0 & \text{if}\;s\in Z(f-g_{0})\cap Z(g).
		\end{cases}
	\end{equation*}
As done earlier, in order to prove that $\phi\in L^{\infty}(\mu,X^{*})$, it is sufficient to show that $\phi$ is strongly measurable. Consider a sequence $\{g_{n}\}$ of simple measurable functions such that $g_{n}(s)\rightarrow g(s)$ for a.e. \enquote{$s$}. For $B_n=\{s\in S: g_n(s)\neq0\}$, define\\
 $\phi_{n}:S\rightarrow X^{*}$ as
	\begin{equation*}
		\phi_{n}(s)=
		\begin{cases}
			F_{\psi_{n}(s)} & \text{if}\;s\in Z(f-g_{0})^c\cap A_n,\\
			cF_{g_{n}(s)} & \text{if}\;s\in Z(f-g_{0})\cap {Z(g)}^c\cap B_n, \\
			0, & elsewhere.
		\end{cases}
	\end{equation*}
	 Clearly, $\{\phi_{n}\}$ is a sequence of simple measurable functions. Consider an $s\in\;S$ such that $\underset{n}\lim \psi_{n}(s)=f(s)-g_{0}(s)$ and $\underset{n}\lim g_{n}(s)=g(s)$, we claim that $\underset{n}\lim \phi_n(s)=\phi(s)$.
	If $s\in Z(f-g_{0})^c$, then $\underset{n}\lim \phi_n(s)=\phi(s)$, as done in Case(i).
	If $s\in Z(f-g_0)\cap Z(g)^c$, then $g(s)\neq0$, so that $s\in B_{n}\;\forall\;n\geq n_0$, for some $n_0\in\;\mathbb{N}$. Thus for $n\geq n_0$, $\phi_n(s)=cF_{g_{n}(s)}$.
	Again, by \Cref{continuity}, the sequence $\{F_{g_{n}(s)}\}_{n\geq n_0}$ converges to $F_{g(s)}$  and hence the sequence $\{\phi_n(s)\}_{n\geq n_0}$ converges to $\phi(s)$.
Lastly, if $s\in Z(f-g_0)\cap Z(g)$, then $\phi_n(s)=0$ for all $n$ and we are done.
	
 Finally, once again, observe that $\int_{S}{\phi(s)(g(s))}\,ds=0$ and as done earlier in Case(i) we have $\|f-g_0\|_1\leq \|f-g_0-g\|_1$.\\
 Hence, $g_{0}\in \mathcal{P}_G(f)$, which completes the proof.
\end{proof}
	       Kripke and Rivlin \cite[Corollary 1.4]{kripke} and Singer \cite[Theorem I.1.7]{singer} established the following:
	       \begin{propn}\label{scalarbjL1}
	        For $f, g\in L^{1}(\mu)$, $f\perp_{BJ}g$ if and only if
	       		       		$$\bigg|\int\limits\limits\limits_{ Z(f)^c}{g(s)\overline{sign}(f(s))} \,ds\bigg|\leq \int\limits\limits\limits_{Z(f)}{|g(s)| \,ds}.$$
	       	\end{propn}
As a direct application of Theorem \ref{L1chr} and Theorem \ref{L1chreal}, we provide its analogue for vector-valued integrable functions.

	       \begin{cor}\label{L1bj}
Let $X$ be a Banach space and $f,g\in L^{1}(\mu,X)$. Then $f\perp_{BJ}g$ in  $L^{1}(\mu,X)$ if and only if 
	       		\begin{equation*}
	       			\hspace{1.5cm}		\bigg|\int\limits\limits\limits_{Z(f)^c}{Re(F_{f(s)}( \alpha g(s))) \,ds}\bigg|\leq\int\limits\limits\limits_{Z(f)} |\alpha| \| g(s)\| \,ds,\ \forall \ \alpha\in\mathbb{C},
	       		\end{equation*} 
       		when $X$ is a complex smooth Banach space, or 
       		     	\begin{equation*}
       		     	\hspace{3cm}	\bigg|\int\limits_{Z(f)^c}{F_{f(s)}(g(s)) \,ds}\bigg|\leq\int\limits\limits\limits_{Z(f)} \| g(s)\| \,ds,
       		     	\end{equation*}
       	     	when $X$ is a real Banach space whose norm is Fr$\acute{e}$chet differentiable. 	
	       \end{cor}
	  We next move on to investigate the elements of best approximation for vector-valued $p$-integrable functions, $1<p<\infty$. For this, we use a result of Leonard \cite[Theorem 3.1]{leo}, where in he proved that for real Banach space $X$, the space $L^{p}(\mu,X)$ is smooth if and only if $X$ is smooth. This result can be proved for complex Banach spaces on the similar lines.
\begin{thm}\label{balp}
	Let $X$ be a Banach space whose norm is Fr$\acute{e}$chet differentiable and let $G$ be a subspace of $L^{p}(\mu,X)$, where $1<p<\infty$. For $f\in L^{p}(\mu,X)\setminus \overline{G},\; g_{0}\in \mathcal{P}_G(f)$ if and only if 
	$$\int\limits\limits\limits_{{Z(f-g_{0})^c}}{\|f(s)-g_{0}(s)\|^{p-1}F_{f(s)-g_{0}(s)}(g(s))\,ds}=0,\;\text{for all} \ \;g\in G.$$
\end{thm}
\begin{proof}
Since the space $L^{p}(\mu,X)$ is smooth, $g_{0}\in \mathcal{P}_G(f)$ if and only if $F_{f-g_{0}}(g)=0$ for every $g\in G$ \cite[Corollary I.1.4]{singer}. So, our primary goal is to determine the  support map at $f-g_{0}$.  For this, let $h\in\;L^p(\mu,X)$ be an arbitrary element. Define a map $\phi_h:S\rightarrow \mathbb{F}$ as 
\begin{equation*}
	\phi_h(s)= 
	\begin{cases}
		\frac{\|f(s)-g_{0}(s)\|^{p-1}}{\|f-g_{0}\|^{p-1}}F_{f(s)-g_{0}(s)}(h(s)) & \text{if}\; f(s)-g_0(s)\neq0\\
		0 & \text{if}\; f(s)-g_0(s)=0.
	\end{cases} 
\end{equation*}
 We first claim that $\phi_h$ is measurable. Let $\{\psi_{n}\}$ and $\{h_n\}$ be the sequences of simple measurable functions such that $ \underset{n}\lim \psi_n(s)=f(s)-g_0(s)$ and $\underset{n}\lim h_n(s)=h(s)$ for a.e. \enquote{$s$}. Set $A_n=\{s:\psi_n(s)\neq0\}$ and define a sequence $\phi_{n}:S\rightarrow \mathbb{F}$ as
\begin{equation*}
	\phi_n(s)= 
	\begin{cases}
		\frac{\|\psi_{n}(s)\|^{p-1}}{\|f-g_{0}\|^{p-1}}F_{\psi_{n}(s)}(h_n(s)) & \; , s\in Z(f-g_0)^c\cap A_n,\\
		0 & \; , \text{otherwise}.
	\end{cases} 
\end{equation*}
 For $s\in S$ such that $\underset{n}\lim \psi_n(s)=f(s)-g_0(s)$ and $\underset{n}\lim  h_n(s)=h(s)$, we claim that $\underset{n}\lim \phi_n(s)=\phi_h(s)$. If $s\in\;Z(f-g_0)^c$, then there exists $n_0\in\mathbb{N}$ such that $s\in A_n$ for all $n\geq n_0.$  Hence, by \Cref{continuity} and using the fact that $p>1$,  the sequence $\{\phi_n(s)\}_{n\geq n_0}$ converges to $\phi_h(s)$. If $s\in\;Z(f-g_0)$, then $|\phi_n(s)|\leq\frac{\|\psi_n(s)\|^{p-1}}{\|f-g_{0}\|^{p-1}}\|h_n(s)\|\;\forall\;n\in\mathbb{N}$, since $\|F_{\psi_n(s)}\|=1$ for $s\in A_n$. Now, $p>1$ implies that $\lim \phi_n(s) = \phi_h(s)=0$. Thus, $\phi_h$ is the a.e. pointwise limit of a sequence of simple measurable functions and by redefining $\phi_{h}$, if required, we derive that $\phi_{h}$ is a measurable function. 

Now define a map $T: L^{p}(\mu,X) \to \mathbb{F}$ as
$$ 	T(h) = \int\limits_{S}{\phi_h(s)\,ds}. $$
To see  that $T$ is well defined consider  $h\in\;L^{p}(\mu,X)$, then
\begin{align*}
	|T(h)|& \leq \frac{1}{\|f-g_{0}\|^{p-1}}\int\limits_{Z(f-g_0)^c}{\|f(s)-g_{0}(s)\|^{p-1}\|h(s)\|\,ds}\\
	& \leq \frac{1}{\|f-g_{0}\|^{p-1}}\bigg(\int\limits_{S}{\|f(s)-g_{0}(s)\|^{p}\,ds}\bigg)^\frac{1}{q}\bigg(\int\limits_{S}{\|h(s)\|^p\,ds}\bigg)^\frac{1}{p}\\
	&=\frac{1}{\|f-g_{0}\|^{p-1}}\|f-g_0\|^\frac{p}{q}\|h\|\\
	&=\|h\|
\end{align*}
 where we have used the fact that the map $s\mapsto \|f(s)-g_0(s)\|^{p-1}$ is in $L^{q}(\mu)$,  $\frac{1}{p}+\frac{1}{q}=1$. Thus, $T\in\;L^{p}(\mu,X)^{*}$ with $\|T\|\leq1$. It is easy to verify that $T(f-g_{0})=\|f-g_{0}\|$ and hence $\|T\|=1$. Since $L^{p}(\mu,X)$ is smooth, $T = F_{f-g_{0}}$ and  this completes the proof.   
\end{proof}
 From \cite[Theorem I.1.11]{singer}, it is known that  for  $f, g\in L^{p}(\mu)$, $1<p<\infty$,  $f\perp_{BJ}g$ if and only if	
 \begin{equation}\label{3}
	\hspace{3cm} \int\limits_{S}{g(s)|f(s)|^{p-1}\overline{\sign}(f(s))} \,ds=0.
		\end{equation}
 Using \Cref{balp}, one can easily obtain its analogue in $L^p(\mu,X)$. 
\begin{cor}\label{Lpbj}
  Let $X$ be a Banach space whose norm is Fr$\acute{e}$chet differentiable and  let $\mu$ be a positive complete measure space. For $f,g\in L^{p}(\mu,X)$, $1<p<\infty$, $f\perp_{BJ}g$ in  $L^{p}(\mu,X)$ if and only if 
			$$\int\limits\limits\limits_{{Z(f)^c}}{\|f(s)\|^{p-1}F_{f(s)}(g(s))\,ds}=0.$$
	\end{cor}

As a consequence, we next deduce an alternate proof of a result of Light \cite[Corollary 2]{light}, which he established using a different technique. 

\begin{thm}\label{subspace}
	Let $f\in L^{p}(\mu,X)$, $(1< p<\infty)$ and let $Y$ be a closed subspace of $X$, where $X$ is a Banach space whose norm is Fr$\acute{e}$chet differentiable and  $\mu$ is a finite complete positive measure. Then $f\bj L^{p}(\mu,Y)$ if and only if $f(s)\bj Y$ for a.e.\enquote{$s$}. 
\end{thm}
\begin{proof}
	First suppose that $f\bj L^{p}(\mu,Y)$. Let, if possible, $f(s)\nbj Y$ for a.e.\enquote{$s$}. Then the set $K=\{s\in S:f(s)\nbj Y\}$ has non-zero measure. Observe that $K=\underset{y\in Y}{\cup}K_{y}$, where $ K_{y}=\{s\in S:f(s)\nbj y\}$. Let $y\in Y$ be such that $K_{y}\neq \emptyset$. We  claim that $K_{y}$ is $\mu$-measurable. Define a function $\phi:S\rightarrow \mathbb{F}$ as
	$$		\phi(s)=
		\begin{cases}
			\|f(s)\|F_{f(s)}(y) & \text{if}\;f(s)\neq 0\\
			0 & \text{if}\; f(s)=0
		\end{cases}
	$$
	Then $\phi$ is $\mu$-measurable. To see this, using the strong measurability of $f$, there exists a sequence $\{f_{n}\}$ of simple measurable functions such that $f_{n}(s)\rightarrow f(s)$ for a.e.\enquote{$s$}. For each $n$, the map $\phi_{n}:S\rightarrow \mathbb{F}$ defined by
	$$		\phi_{n}(s)=
		\begin{cases}
			\|f_{n}(s)\|F_{f_{n}(s)}(y) & \text{if}\;f_{n}(s)\neq 0\\
			0 & \text{if}\; f_{n}(s)=0
		\end{cases}
	$$
	is a simple measurable function.
	 Consider $s\in S$ for which $\underset{n \to \infty}\lim f_{n}(s) = f(s)$.
	
	If $f(s)\neq0$, since the norm is Fr$\acute{e}$chet differentiable, by \Cref{continuity} it is easy to see $\underset{n \to \infty}\lim F_{f_{n}(s)} = F_{f(s)}$ (take a tail of the sequence $\{f_{n}\}$, if required) in the operator norm topology. Thus, $\underset{n \to \infty}\lim F_{f_{n}(s)}(y)= F_{f(s)}(y)$ which gives $\underset{n \to \infty}\lim \phi_{n}(s)=\phi(s)$.
	
	If $f(s)=0$,  then using the fact that $\|F_{f_{n}(s)}\|=1$, we have $|\phi_{n}(t)|\leq\|f_{n}(t)\|$ for all $t \in S$ and  $n \in \mathbb{N}$. 
	 This gives that $\underset{n \to \infty}\lim \phi_{n}(s) =0 $. Thus, in both the cases, the function $\phi$ is a.e. limit of a sequence of simple measurable functions and hence $\phi$ is measurable being the measure is complete.
	  By James criteria and the smoothness of $X$, $K_{y}=\{s\in S:F_{f(s)}(y)\neq0\}$. Further, by using the fact that $f(s)\neq0$ for all $s\in K_{y}$, we observe that $K_{y}=\phi^{-1}(\mathbb{F}\setminus\{0\})$ and hence measurable.\\ 
	  {\bf Case(i): $X$ is a real Banach space.} Fix a $y\in Y$ such that $\mu(K_{y})\neq0$ and write $K_{y}=K_{y}^{+}\cup K_{y}^{-}$ where $K_{y}^{+}=\{s\in S:F_{f(s)}(y)>0\}$ and $K_{y}^{-}= \{s\in S:F_{f(s)}(y)<0\}$. Again, observe that both $K_{y}^{+}$ and $K_{y}^{-}$ are measurable as $K_{y}^{+}=\phi^{-1}(0,\infty)$ and $K_{y}^{-}=\phi^{-1}(-\infty,0)$.
	  Without loss of generality, assume that $\mu(K_{y}^{+})\neq0$ and define a map $g:S\rightarrow Y$ as $g(s)=y\chi_{K_{y}^{+}}(s)$, where $\chi_{K_{y}^{+}}$ denotes the characteristic function on $S$.
		Then $g\in L^{p}(\mu,Y)$, as $\mu$ is a finite measure. Also
	$$\int\limits_{{Z(f)}^c}{\|f(s)\|^{p-1}F_{f(s)}(g(s))\,ds} =\int\limits\limits\limits_{K_{y}^{+}}{\|f(s)\|^{p-1}F_{f(s)}(y)\,ds} \neq 0	$$
	 which, by \Cref{Lpbj}, is a contradiction to the fact that $f\bj g$.\\
		{\bf  Case(ii): $X$ is a complex Banach space.} Fix a $y\in\;Y$ such that $\mu(K_y)\neq0$. Since $F_{f(s)}(y)=(ReF_{f(s)})(y)-i(ReF_{f(s)})(iy)$, therefore,
		 \begin{align*}
		  K_{y}&=\{s\in S: (ReF_{f(s)})(y)-i(ReF_{f(s)})(iy)\neq0\}\\
		  &=K_{1}\cup K_{2}\cup K_{3}\cup K_{4}
		  \end{align*}
	  where, $K_{1}=\{s\in S: (ReF_{f(s)})(y)\geq0\},\;K_{2}=\{s\in S: (ReF_{f(s)})(y)\leq0\}$ and $K_{3}=\{s\in S: (ReF_{f(s)})(iy)\geq0\},\;K_{4}=\{s\in S: (ReF_{f(s)})(iy)\leq0\}$.  First, observe that $K_{1}$ is a measurable set. For this, define  $\psi:S\rightarrow \mathbb{R}$ as
	  $$		\psi(s)=
	  \begin{cases}
	  	\|f(s)\|(ReF_{f(s)})(y) & \text{if}\;f(s)\neq 0\\
	  	0 & \text{if}\; f(s)=0.
	  \end{cases}
	  $$ 
	  As done earlier, the function $\psi$ is the a.e. limit of the sequence of simple measurable functions and hence $\psi$ is a measurable. Thus, $K_{1}=\psi^{-1}((0,\infty))$ is a measurable set. Similarly, one can verify that each $K_{i},\;i\in\{2,3,4\}$ is a measurable set.  Since $\mu(K_{y})\neq0$, without loss of generality, assume $\mu(K_{1})\neq0$ and then proceeding in the same manner as in Case(i), we obtain the desired result.
	  
		Conversely, if $f(s)\bj Y$ for a.e.\enquote{$s$}, then $F_{f(s)}(y)=0$ for a.e. \enquote{$s$} and for each $y\in Y$. Therefore $\int_{S}{\|f(s)\|^{p-1}F_{f(s)}(g(s))\,ds}=0$, for every $g\in L^{p}(\mu,Y)$ and by \Cref{Lpbj}, $f\bj L^{p}(\mu, Y)$. 
	\end{proof}

\section{Birkhoff-James orthogonality and tensor product}
With all the ingredients prepared, we are now ready to discuss BJ-orthogonality in the tensor product spaces $L^{p}(\mu)\otimes^{\Delta_{p}} X$, $1\leq p<\infty$.
\begin{thm}\label{p-norm-result}
	Let $X$ be a Banach space whose norm is Fr$\acute{e}$chet differentiable and $1<p<\infty$. Then $f\otimes x\perp_{BJ}g\otimes y$ in $L^{p}(\mu)\otimes^{\Delta_{p}}X$ if and only if either $f\perp_{BJ}g$ in $L^{p}(\mu)$ or $x\perp_{BJ}y$ in $X$.
\end{thm}
\begin{proof}
 Let $f\otimes x\perp_{BJ}g\otimes y$. We assume $x$ to be non-zero. Since $L^{p}(\mu)\otimes^{\Delta_p}X$ is isometrically isomorphic to $L^{p}(\mu,X)$, we have that  $f_x\perp_{BJ} g_y$ in $L^{p}(\mu, X)$ where $f_x$ and $g_y$ correspond to $f\otimes x$ and $g\otimes y$, respectively. Thus, by \Cref{Lpbj}, we have 
	$$			\int\limits\limits\limits_{{Z(f_x)}^c}{\|f_x(s)\|^{p-1}F_{f_x(s)}(g_{y}(s))\,ds}=0,$$	
		which further implies
		$$				 \int\limits\limits\limits_{{Z(f_x)}^c}{|f(s)|^{p-1}\|x\|^{p-1}F_{f(s)x}(g(s)y)\,ds}=0.$$	
	  Since $\frac{\overline{f(s)}}{|f(s)|}F_{x}(f(s)x)=\|f(s)x\|$ for $0\neq f(s)$ 
	and the space $X$ is smooth, we have $F_{f(s)x}=\frac{\overline{f(s)}}{|f(s)|}F_{x}$ for $0\neq f(s)$. Observing that $Z(f_x) = Z(f)$, the above equation becomes
			$$
			 (F_{x}(y)\|x\|^{p-1})\int\limits\limits\limits_{{Z(f)}^c}{|f(s)|^{p-1}g(s)\overline{\sign(f(s))}\,ds}=0.
		$$
	Thus, either $\int\limits\limits\limits_{{Z(f)}^c}{|f(s)|^{p-1}g(s)\overline{\sign(f(s))}\,ds}=0$ or  $F_{x}(y)\|x\|^{p-1}=0$. Since $X$ is smooth, by \Cref{3}, either  $f\perp_{BJ}g$ or $x\perp_{BJ}y$. Converse follows from \cite[Theorem 3.1]{jain}.
\end{proof}
 It is interesting to note that in Theorem \ref{p-norm-result} if we take $X$ to be a complex Banach space and $p=1$, the conclusion may not hold as seen in the following example.

\begin{example}\label{example1}
	Consider the measure space $(\mathbb{N}, P(\mathbb{N}), \mu),$ where $\mu$ is the counting measure, $P(\mathbb{N})$ denotes the power set of $\mathbb{N}$, and let $X=\ell^{2}(\mathbb{C}).$ Take $A=\{1,2,3\},B=\{2,3,5\}$ and let $x=(i,-i,0,0,.....)$, $y=(i,0,-i,0,0,...)\in X.$ Now,
	\begin{align*}
	 \int\limits\limits\limits_{{Z(\chi_{A})}^c}{\chi_{B}(s)\overline{\sign(\chi_{A}(s))}\,ds}& =\mu(A\cap B)\\
	 	& >\mu(A^c\cap B)\\
	 	& =\int\limits\limits\limits_{Z(\chi_{A})}{|\chi_{B}(s)|\,ds}. 
	 \end{align*}
	 Thus, by \Cref{scalarbjL1}, $\chi_{A}\nbj\chi_{B}$. Also, $x\nbj y$ since BJ-orthogonality coincides with the usual orthogonality in Hilbert spaces. Now, we claim that $\chi_{A}\ot x\perp_{BJ}\chi_{B}\ot y$ in $L^{1}(\mu)\ot^{\gamma}X$. As done earlier,  it is sufficient to prove that $h_{1}\perp_{BJ}h_{2}$ in $L^{1}(\mu,X)$, where $h_{1},\;h_{2}\in L^{1}(\mu,X)$ correspond to $\chi_{A}\ot x$ and $\chi_{B}\ot y$, respectively. To see this, for any scalar $\alpha\in\mathbb{C}$, we have
	\begin{align*}
		\hspace{-0.9cm}	\bigg| \int\limits\limits\limits_{{Z(h_{1})}^c}Re(F_{h_{1}(s)}(\alpha h_{2}(s)))\,ds\bigg|&=\bigg|\int\limits\limits\limits_{{Z(h_{1})}^c}Re\bigg(\left\langle \alpha h_{2}(s),\frac{h_{1}(s)}{\|h_{1}(s)\|}\right\rangle_{X}\bigg)\,ds\bigg|\\
&=\bigg|\int\limits\limits\limits_{{Z(\chi_{A})}^c}Re\bigg(\left\langle \alpha\chi_{B}(s)y,\frac{\chi_{A}(s)x}{|\chi_{A}(s)|\|x\|}\right\rangle_{X}\bigg)\,ds\bigg|\\
	&=\bigg| \int\limits\limits\limits_{A\cap B}\frac{Re(\alpha(\langle y,\;x \rangle_{X}))}{\|x\|_{X}}\,ds\bigg|\\
	&=\frac{|Re(\alpha \ \mu(A\cap B))|}{\sqrt{2}}\\
	&=\frac{|Re(2\alpha)|}{\sqrt{2}}.
	\end{align*}
	 On the other hand,
		 \begin{align*}
	  \int\limits\limits\limits_{Z(h_{1})}{\|\alpha h_{2}(s)\|\,ds}&=\int\limits\limits\limits_{Z(\chi_{A})}{\|\alpha \chi_{B}(s)y\|_{X}\,ds}\\&=\int\limits\limits\limits_{A^{c}\cap B}{|\alpha|\|y\|_{X}\,ds}\\&=\sqrt{2}|\alpha|\mu(A^{c}\cap B)\\&=\sqrt{2}|\alpha|\\&\geq\bigg| \int\limits\limits\limits_{{Z(h_{1})}^c}Re(F_{h_{1}(s)}(\alpha h_{2}(s)))\,ds\bigg|
	  \end{align*}
	  and hence by \Cref{L1bj}, $h_{1}\perp_{BJ}h_{2}$. This proves the claim.
\end{example}
 Lastly, we present an example to show that $f_{1}\otimes f_{2}\perp_{BJ}g_{1}\otimes g_{2}$ in $L^{1}(\mu)\otimes^{\gamma}L^{1}(\nu)$ need not imply either $f_{1}\perp_{BJ}g_{1}$ or $f_{2}\perp_{BJ}g_{2}.$  
\begin{example}
	Consider the measure spaces $(\mathbb{N}, P(\mathbb{N}), \mu)$ and $(\mathbb{R}, \mathcal{M}(\mathbb{R}), \nu)$, where $\mu$ is the counting measure, $\nu$ is the Lebesgue measure, $P(\mathbb{N})$ denotes the power set of $\mathbb{N}$ and $ \mathcal{M}(\mathbb{R})$ denotes the algebra of Lebesgue measurable subsets of $\mathbb{R}$. Take $A=\{1,2,3\}, B=\{2,3,5\}$ and $C=[-1,2], D=[-2,1]$. As done in \Cref{example1}, since $\mu({A\cap B})>\mu({A^{c}\cap B})$ and $\nu({C\cap D})>\nu({C^{c}\cap D})$, therefore, by \Cref{scalarbjL1} neither $\chi_{A}\perp_{BJ}\chi_{B}$ nor $\chi_{C}\perp_{BJ}\chi_{D}$. Now, we claim that $\chi_{A}\ot \chi_{C}\perp_{BJ}\chi_{B}\ot\chi_{D}$ in $L^{1}(\mu)\otimes^{\gamma} L^{1}(\nu)$. As done earlier, it is sufficient to prove that $\chi_{A}\chi_{C}\perp_{BJ}\chi_{B}\chi_{D}$ in $L^{1}(\mu\times \nu).$ To see this, consider
	$$
		 \hspace*{-7cm} \bigg|\int\limits\limits\limits\limits_{{Z(\chi_{A}\chi_{C})}^c}{\chi_{B}\chi_{D}\overline{\sign}(\chi_{A}\chi_{C}) \,d(\mu\times\nu)\bigg|} $$  
			$$\hspace*{3cm} =\bigg|\bigg(\int\limits\limits\limits_{{Z(\chi_{A})}^c}{\chi_{B}(s)\overline{\sign{\chi_{A}(s)}}\,ds}\bigg)\bigg(\int\limits\limits\limits_{{Z(\chi_{C})}^c}{\chi_{D}(t)\overline{\sign{\chi_{C}(t)}}\,dt}\bigg)\bigg|$$
		$$\hspace{-4.5cm}	=\mu(A\cap B)\nu(C\cap D)$$
	$$\hspace*{-7.5cm}=4.$$
	On the other hand
	$$\hspace*{-10cm} \int\limits_{Z(\chi_{A}\chi_{C})}{|\chi_{B}\chi_{D}| \, d(\mu\times\nu)}$$
	$$=\int\limits\limits\limits_{C^c}\int\limits\limits\limits_{A}{|\chi_{B}(s)\chi_{D}(t)| \, ds\,dt}+\int\limits\limits\limits_{C}\int\limits\limits\limits_{A^c }{|\chi_{B}(s)\chi_{D}(t)| \, ds\,dt}
	+\int\limits\limits\limits_{C^c}\int\limits\limits\limits_{A^c}{|\chi_{B}(s)\chi_{D}(t)| \, ds\,dt}$$
	$$\hspace{-.75cm} =\mu(A\cap B)\nu(C^{c}\cap D)+\mu(A^{c}\cap B)\nu(C\cap D)
+\mu(A^{c}\cap B)\nu(C^{c}\cap D) = 5.$$
	
		
 Thus, by \Cref{scalarbjL1},	 $\chi_{A}\chi_{C}\perp_{BJ}\chi_{B}\chi_{D},$ which proves the claim.
	\end{example}\vspace{0.5cm}


\begin{thebibliography}{00}
	\bibitem{birk} Birkhoff, G. \emph{Orthogonality in linear metric spaces}, Duke Math. J. 1 (1935), 169-172.
	\bibitem{cudia} Cudia, D.~F. \emph{The geometry of Banach spaces, smoothness}, Trans. Amer. Math. Soc. 110 (1964), 284-314.
	\bibitem{floret} Defant, A. and Floret, K. \emph{Tensor norms and operator ideals}, Elsevier Science Publishers B.V., 1993.
	\bibitem{giles} Giles, J.R. \emph{On a characterization of differentiability of the norm of a normed linear space}, Journal of the Australian Mathematical Society 12(1) (1971), 106-114.
	\bibitem{gro} Grothendieck A. {\emph Resum$\acute{e}$ de la th$\acute{e}$orie m$\acute{e}$trique des produits tensorielles topologiques}, Bol. Soc. Math. Sao. Paolo 8(1953), 1-79.
	\bibitem{james} James, R.~C. \emph{Orthogonality and linear functionals in normed linear spaces,} Trans. Amer. Math. Soc. 61 (1947), 265-292.
	
	\bibitem{keckic} Ke$\check{c}$ki$\grave{c}$, D.~J. \emph{G$\hat{a}$teaux derivative of $B(H)$ norm,} Proc. Amer. Math. Soc. 133 (2005), 2061-2067.
	
	\bibitem{kripke} Kripke, B.~R. and Rivlin, T.~J. \emph{Approximation in the metric of $L^1(X,\mu)$,} Trans. Amer. Math. Soc. 119 (1965), 101-122.
	
	\bibitem{leo} Leonard, I.~E. \emph{ Smoothness and Duality in $L^{p}(E,\mu)$}, J. Math. Anal. Appl. 46 (1974),  513-522.
	
\bibitem{light} Light, W.~A.  \emph{Proximinality in $L^{p}(\mu,Y)$}, Rocky Mountain J. Math. 19 (1989), 251-259.

	\bibitem{jain} Mohit and Jain, R. \emph{Birkhoff-James orthogonality in certain tensor products of Banach spaces}, Operators and Matrices 17(1) (2023), 235-244.
	\bibitem{singer} Singer, I. \emph{Best approximation in normed linear spaces by elements of linear subspaces}, Springer-Verlag, Berlin, 1970.

	\bibitem{smirnov} Smirnov, G.~S. \emph{A remark on the best approximation in the mean of vector valued functions,} Ukr. Math. J.  (1989), 703-704.
	

	\end{thebibliography}
\end{document}